\documentclass{amsproc}%
\usepackage{graphicx}
\usepackage{amscd}
\usepackage{amsmath}
\usepackage{amsfonts}
\usepackage{amssymb}%
\theoremstyle{plain}

\numberwithin{equation}{section}

\begin{document}
\title[Embedding truncated skew polynomial rings into matrix rings]{Embedding truncated skew polynomial rings into matrix rings and embedding of a
ring into $2\times2$\ supermatrices}
\author{Jen\H{o} Szigeti}
\address{Institute of Mathematics, University of Miskolc, Miskolc, Hungary 3515}
\email{matjeno@uni-miskolc.hu}
\thanks{The author was supported by OTKA K-101515 of Hungary}
\thanks{This research was carried out as part of the TAMOP-4.2.1.B-10/2/KONV-2010-0001
project with support by the European Union, co-financed by the European Social Fund.}

\begin{abstract}
For an endomorphism $\sigma:R\longrightarrow R$ with $\sigma^{t}=1$ we prove
that the truncated polynomial ring (algebra) $R[w,\sigma]/(w^{t})$ embeds into
$\mathrm{M}_{t}(R[z]/(z^{t}))$. For an involution $\sigma$ we exhibit an
embedding $R\longrightarrow\mathrm{M}_{2,1}^{\sigma}(R)$, where $\mathrm{M}%
_{2,1}^{\sigma}(R)$ is the algebra of the so called $(\sigma,2,1)$\ supermatrices.

\end{abstract}
\subjclass{16S36,16S50,16R10,16R20}
\keywords{skew polynomial algebra, truncated skew polynomial algebra, matrix algebra,
embedding, supermatrix defined by an involution }
\maketitle

\noindent1. INTRODUCTION

\bigskip

The main inspiration of the present work are the various embedding results in
[SvW] and [MMSvW]. A certain embedding of the two-generated Grassmann algebra
$E^{(2)}$ into a $2\times2$ matrix algebra over a commutative ring leads to a
Cayley-Hamilton identity of degree $2n$ for any $n\times n$ matrix over
$E^{(2)}$ (see [SvW]). For a field $K$ (of characteristic zero) let%
\[
E^{(m)}=K\left\langle v_{1},\ldots,v_{m}\mid v_{i}v_{j}+v_{j}v_{i}=0\text{ for
all }1\leq i\leq j\leq m\right\rangle
\]
denote the $m$-generated Grassmann algebra and $\mathrm{M}_{n}(R)$ denote the
full $n\times n$ matrix ring ($K$-algebra) over a ring ($K$-algebra) $R$ with
identity $I_{n}\in\mathrm{M}_{n}(R)$. In [MMSvW] a so called constant trace
representation ($K$-embedding)%
\[
\varepsilon^{(m)}:E^{(m)}\longrightarrow\mathrm{M}_{2^{m-1}}(K[z_{1}%
,\ldots,z_{m}]/(z_{1}^{2},\ldots,z_{m}^{2})),
\]
is presented, where the ideal $(z_{1}^{2},\ldots,z_{m}^{2})$ of the
commutative polynomial algebra $K[z_{1},\ldots,z_{m}]$ is generated by the
monomials $z_{1}^{2},\ldots,z_{m}^{2}$. One of the remarkable consequences of
this embedding is a Cayley-Hamilton identity (with coefficients in $K$) of
degree $2^{m-1}n$ for $n\times n$ matrices over $E^{(m)}$.

The induction step in the construction of the above $\varepsilon^{(m)}$ is
based on the observation that $E^{(m)}[w,\tau]/(w^{2})\cong E^{(m+1)}$ as
$K$-algebras, where $\tau:E^{(m)}\longrightarrow E^{(m)}$ is the natural
involution defined by the well known $\mathbb{Z}_{2}$-grading $E^{(m)}%
=E_{0}^{(m)}\oplus E_{1}^{(m)}$ and $(w^{2})$ is the ideal of the skew
polynomial algebra $E^{(m)}[w,\tau]$ generated by $w^{2}$. The main ingredient
of the mentioned induction is a "Fundamental Embedding" $\mu:R[w,\sigma
]/(w^{2})\longrightarrow\mathrm{M}_{2}(R[z]/(z^{2}))$, which is defined for an
arbitrary involution $\sigma:R\longrightarrow R$. In Section 2 we give a far
reaching generalization of this "Fundamental Embedding". We note that the idea
of considering the truncated polynomial ring ($K$-algebra) $R[w,\sigma
]/(w^{2})$ comes from [SSz].

The use of an endomorphism (involution) $\sigma:R\longrightarrow R$ enables us
to give a generalization of the concept of a supermatrix. Supermatrices over
the infinitely generated Grassmann algebra play an important role in Kemer's
classification of T-prime T-ideals (see [K]). Section 3 contains an embedding
of $R$ into a $2\times2$\ supermatrix algebra (over $R$) determined by
$\sigma$.

\bigskip

\noindent2. TRUNCATED SKEW\ POLYNOMIAL\ RINGS AND EMBEDDINGS

\bigskip

For a ring ($K$-algebra) endomorphism $\sigma:R\longrightarrow R$\ let us
consider the skew polynomial ring ($K$-algebra) $R[w,\sigma]$ in the skew
indeterminate $w$. The elements of $R[w,\sigma]$ are left polynomials of the
form $f(w)=r_{0}+r_{1}w+\cdots+r_{k}w^{k}$ with $r_{0},r_{1},\ldots,r_{k}\in
R$. Besides the obvious addition, we have the following multiplication rule in
$R[w,\sigma]$: $wr=\sigma(r)w$ for all $r\in R$ and
\[
(r_{0}+r_{1}w+\cdots+r_{k}w^{k})(s_{0}+s_{1}w+\cdots+s_{l}w^{l})=u_{0}%
+u_{1}w+\cdots+u_{k+l}w^{k+l},
\]
where%
\[
u_{m}=\underset{i+j=m,i\geq0,j\geq0}{\sum}r_{i}\sigma^{j}(s_{j})\text{
}(\text{for }0\leq m\leq k+l).
\]

If $\sigma^{t}=1$ (such a $\sigma$\ is an automorphism), then $w^{t}$ is a
central element of $R[w,\sigma]$: we have $\sigma^{t}(r)=r$ and $w^{t}%
r=w^{t-1}\sigma(r)w=\cdots=\sigma^{t}(r)w^{t}=rw^{t}$ for all $r\in R$,
moreover $w^{t}$ commutes with the powers of $w$. Thus the ideal $(w^{t})\lhd
R[w,\sigma]$ generated by $w^{t}$ can be written as $(w^{t})=R[w,\sigma
]w^{t}=w^{t}R[w,\sigma]$. For any element $f(w)+(w^{t})$ of the truncated
polynomial ring ($K$-algebra) $R[w,\sigma]/(w^{t})$ there exists a unique
sequence of coefficients $r_{0},r_{1},\ldots,r_{t-1}\in R$ such that%
\[
r_{0}+r_{1}w+\cdots+r_{t-1}w^{t-1}+(w^{t})=f(w)+(w^{t}).
\]
Hence the elements of $R[w,\sigma]/(w^{t})$ can be represented by left
polynomials of degree less or equal than $t-1$.

The following is called "Fundamental Embedding" in [MMSzvW].

\bigskip

\noindent\textbf{2.1.Theorem} ([MMSzvW])\textbf{.}\textit{ For an involution
}$\sigma:R\longrightarrow R$\textit{, putting}%
\[
\mu(r_{0}+r_{1}w+(w^{2}))=\left[
\begin{array}
[c]{cc}%
r_{0}+(z^{2}) & r_{1}z+(z^{2})\\
\sigma(r_{1})z+(z^{2}) & \sigma(r_{0})+(z^{2})
\end{array}
\right]
\]
\textit{(with} $r_{0},r_{1}\in R$\textit{) gives an embedding }$\mu
:R[w,\sigma]/(w^{2})\longrightarrow\mathrm{M}_{2}(R[z]/(z^{2}))$\textit{.}

\bigskip

Now we present the following generalization of Theorem 2.1 (already announced
in [MMSzvW]).

\bigskip

\noindent\textbf{2.2.Theorem.}\textit{ For an endomorphism }$\sigma
:R\longrightarrow R$\textit{ with }$\sigma^{t}=1$\textit{,} \textit{putting}%
\[
\mu(r_{0}+r_{1}w+\cdots+r_{t-1}w^{t-1}+(w^{t}))=\left[  \sigma^{i-1}%
(r_{j-i})z^{j-i}+(z^{t})\right]  _{t\times t}%
\]
\textit{gives an embedding }$\mu:R[w,\sigma]/(w^{t})\longrightarrow
\mathrm{M}_{t}(R[z]/(z^{t}))$\textit{,} \textit{where the difference }
$j-~i\in\{0,1,\ldots,t-1\}$\textit{ is taken modulo }$t$\textit{, and the
element of the factor algebra} $R[z]/(z^{t})$\textit{ in the }$(i,j)$\textit{
position of the }$t\times t$\textit{ matrix }$\left[  \sigma^{i-1}%
(r_{j-i})z^{j-i}\!+\!(z^{t})\right]  _{t\times t}$ \textit{is} $\sigma
^{i-1}(r_{j-i})z^{j-i}\!+\!~(z^{t})$\textit{. The trace of }$\left[
\sigma^{i-1}(r_{j-i})z^{j-i}\!+\!(z^{t})\right]  _{t\times t}$ \textit{is in
the fixed ring }$R^{\sigma}\!=\!\{r\in R\!\mid\!\sigma(r)=r\}$\textit{:}%
\[
\mathrm{tr}(\mu(r_{0}+r_{1}w+\cdots+r_{t-1}w^{t-1}+(w^{t})))=\mathrm{tr}%
(\left[  \sigma^{i-1}(r_{j-i})z^{j-i}+(z^{t})\right]  _{t\times t})=
\]%
\[
r_{0}+\sigma(r_{0})+\cdots+\sigma^{t-1}(r_{0})+(z^{t})\in R^{\sigma}+(z^{t}).
\]

\bigskip

\noindent\textbf{Proof.} We only have to prove the multiplicative property of
$\mu$.

\noindent In order to avoid confusion, for $i,j\in\{1,\ldots,t\}$ let
$j\circleddash i$ denote the modulo $t$ difference:%
\[
j\circleddash i=\left\{
\begin{array}
[c]{c}%
j-i\text{ if }i\leq j\\
(j-i)+t\text{ if }j\leq i-1
\end{array}
\right.  .
\]
The $(p,q)$ entry in the product of the $t\times t$ matrices%
\[
\mu(r_{0}+r_{1}w+\cdots+r_{t-1}w^{t-1}+(w^{t}))=\left[  \sigma^{i-1}%
(r_{j\circleddash i})z^{j\circleddash i}+(z^{t})\right]  _{t\times t}%
\]
and%
\[
\mu(s_{0}+s_{1}w+\cdots+s_{t-1}w^{t-1}+(w^{t}))=\left[  \sigma^{i-1}%
(s_{j\circleddash i})z^{j\circleddash i}+(z^{t})\right]  _{t\times t}%
\]
is%
\[
a_{p,q}=\underset{j=1}{\overset{t}{%
{\textstyle\sum}
}}\left(  \sigma^{p-1}(r_{j\circleddash p})z^{j\circleddash p}+(z^{t})\right)
\left(  \sigma^{j-1}(s_{q\circleddash j})z^{q\circleddash j}+(z^{t})\right)  =
\]%
\[
\left(  \underset{j=1}{\overset{t}{%
{\textstyle\sum}
}}\sigma^{p-1}(r_{j\circleddash p})\sigma^{j-1}(s_{q\circleddash j})\right)
z^{(j\circleddash p)+(q\circleddash j)}+(z^{t}).
\]
Since $\sigma^{t}=1$, we have $\sigma^{p-1}(\sigma^{j\circleddash
p}(s_{q\circleddash j}))=\sigma^{j-1}(s_{q\circleddash j})$. It is
straightforward to check that if $(j\circleddash p)+(q\circleddash j)\leq t-1$
holds, then%
\[
(j\circleddash p)+(q\circleddash j)=q\circleddash p\text{ and }0\leq
j\circleddash p\leq q\circleddash p.
\]
In view of the above observations and using $i=j\circleddash p$, we can see
that%
\[
a_{p,q}=\sigma^{p-1}\left(  \underset{j=1}{\overset{t}{%
{\textstyle\sum}
}}r_{j\circleddash p}\sigma^{j\circleddash p}(s_{q\circleddash j})\right)
z^{(j\circleddash p)+(q\circleddash j)}+(z^{t})=
\]%
\[
\sigma^{p-1}\left(  \underset{i=0}{\overset{q\circleddash p}{%
{\textstyle\sum}
}}r_{i}\sigma^{i}(s_{(q\circleddash p)-i})\right)  z^{q\circleddash p}%
+(z^{t})=\sigma^{p-1}(u_{q\circleddash p})z^{q\circleddash p}+(z^{t})
\]
is the $(p,q)$ entry of $\mu(u_{0}+u_{1}w+\cdots+u_{t-1}w^{t-1}+(w^{t}))$,
where%
\[
u_{0}+u_{1}w+\cdots+u_{t-1}w^{t-1}+(w^{t})=
\]%
\[
(r_{0}+r_{1}w+\cdots+r_{t-1}w^{t-1}+(w^{t}))(s_{0}+s_{1}w+\cdots
+s_{t-1}w^{t-1}+(w^{t})).
\]
holds in $R[w,\sigma]/(w^{t})$.

Since%
\[
\sigma(r_{0}+\sigma(r_{0})+\cdots+\sigma^{t-1}(r_{0}))=\sigma(r_{0}%
)+\sigma^{2}(r_{0})+\cdots+\sigma^{t-1}(r_{0})+\sigma^{t}(r_{0})
\]
and $\sigma^{t}(r_{0})=r_{0}$, the trace of $\left[  \sigma^{i-1}%
(r_{j-i})z^{j-i}+(z^{t})\right]  _{t\times t}$ is in $R^{\sigma}+(z^{t})$.
$\square$

\bigskip

Now consider the free associative $K$-algebra $K\left\langle x_{1}%
,\ldots,x_{m},\ldots\right\rangle $ generated by the (non-commuting)
indeterminates $x_{1},\ldots,x_{m},\ldots$and let
\[
S_{m}(x_{1},\ldots,x_{m})=\sum_{\pi\in\mathrm{Sym}(m)}\mathrm{sgn}(\pi
)x_{\pi(1)}\cdots x_{\pi(m)}%
\]
be the standard polynomial in $K\left\langle x_{1},\ldots,x_{m},\ldots
\right\rangle $. In the following corollaries we keep the notations and the
conditions of Theorem 2.2.

\bigskip

\noindent\textbf{2.3.Corollary.}\textit{ If }$R$\textit{ is commutative and
}$n\geq1$\textit{ is an integer, then }$S_{2tn}=0$\textit{ is an identity on
}$\mathrm{M}_{n}(R[w,\sigma]/(w^{t}))$\textit{. In particular }$S_{2t}%
=0$\textit{ is an identity on }$R[w,\sigma]/(w^{t})$\textit{.}

\bigskip

\noindent\textbf{Proof.} The natural extension%
\[
\mu_{n}:\mathrm{M}_{n}(R[w,\sigma]/(w^{t}))\longrightarrow\mathrm{M}%
_{n}(\mathrm{M}_{t}(R[z]/(z^{t})))\cong\mathrm{M}_{tn}(R[z]/(z^{t}))
\]
of $\mu$ is an embedding. Since $R[z]/(z^{t})$ is also commutative and
$S_{2tn}=0$ is an identity on $\mathrm{M}_{tn}(R[z]/(z^{t}))$ by the
Amitsur-Levitzki theorem (see [Dr, DrF]), the proof is complete. $\square$

\bigskip

\noindent\textbf{2.4.Corollary.}\textit{ If }$S_{m}=0$\textit{ is an identity
on }$R$\textit{ and }$n\geq1$\textit{ is an integer, then }$S_{(m-1)t^{2}%
n^{2}+1}=0$\textit{ is an identity on }$\mathrm{M}_{n}(R[w,\sigma]/(w^{t}%
))$\textit{. In particular }$S_{(m-1)t^{2}+1}=0$\textit{ is an identity on
}$R[w,\sigma]/(w^{t})$\textit{.}

\bigskip

\noindent\textbf{Proof.} Since $S_{m}=0$ is also an identity on $R[z]/(z^{t}%
)$, using $\mu_{n}$ and Theorem 5.5 of Domokos [Do] completes the proof.
$\square$

\bigskip

\noindent\textbf{2.5.Corollary.}\textit{ If }$R$\textit{ is a PI-algebra and
}$n\geq1$\textit{ is an integer, then }$\mathrm{M}_{n}(R[w,\sigma]/(w^{t}%
))$\textit{ is also a PI-algebra. In particular }$R[w,\sigma]/(w^{t})$\textit{
is a PI-algebra.}

\bigskip

\noindent\textbf{Proof.} Since $R[z]/(z^{t})$ is also PI, using $\mu_{n}$ and
the well known fact that full matrix algebras over a PI-algebra are PI (a
special case of Regev's tensor product theorem), the proof is complete.
$\square$

\bigskip

\noindent\textbf{2.6.Theorem.}\textit{ If }$R$\textit{\ is commutative and
}$\sigma:R\longrightarrow R$\textit{ is an endomorphism, then the fixed ring
(algebra) }$R^{\sigma}=\{r\in R\mid\sigma(r)=r\}$\textit{\ of }$\sigma
$\textit{ is a central subring (algebra) of }$R[w,\sigma]$\textit{ and
}$R^{\sigma}+(z^{t})\subseteq\mathrm{Z}(R[w,\sigma]/(w^{t}))$\textit{. If
}$\sigma^{t}=1$ \textit{and }$A\in\mathrm{M}_{n}(R[w,\sigma]/(w^{t}))$\textit{
is an }$n\times n$\textit{ matrix, then }$A$\textit{ satisfies a
"Cayley-Hamilton" identity of the form}%
\[
A^{tn}+c_{1}A^{tn-1}+\cdots+c_{tn-1}A+c_{tn}I_{n}=0,
\]
\textit{where }$c_{i}\in R^{\sigma}=\{r\in R\mid\sigma(r)=r\}$\textit{,
}$1\leq i\leq tn$\textit{. In particular }$R[w,\sigma]/(w^{t})$\textit{ is
integral over }$R^{\sigma}$\textit{ of degree }$t$\textit{.}

\bigskip

\noindent\textbf{Proof.} The containments $R^{\sigma}\subseteq\mathrm{Z}%
(R[w,\sigma])$ and $R^{\sigma}+(z^{t})\subseteq\mathrm{Z}(R[w,\sigma
]/(w^{t}))$ are clear. In the rest of the proof we follow the steps of the
proof of Theorem 2.1 in [MMSzvW]. Let $A=[a_{i,j}]$ and take $\mu_{n}$ from
the proof of Corollary 2.3. The trace of the $tn\times tn$ matrix $B=\mu
_{n}(A)\in\mathrm{M}_{tn}(R[z]/(z^{t}))$ is the sum of the traces of the
diagonal $t\times t$ blocks:%
\[
\mathrm{tr}(B)=\overset{n}{\underset{i=1}{\sum}}\mathrm{tr}(\mu(a_{i,i})).
\]
Theorem 2.2 ensures that $\mathrm{tr}(\mu(a_{i,i}))\in R^{\sigma}+(z^{t})$ for
each $1\leq i\leq n$. For the sake of simplicity we can take $\mathrm{tr}%
(\mu(a_{i,i}))\in R^{\sigma}$. It follows that $\mathrm{tr}(B)\in R^{\sigma}$.
The coefficients of the characteristic polynomial%
\[
\det(xI-B)=c_{0}x^{tn}+c_{1}x^{tn-1}+\cdots+c_{tn-1}x+c_{tn}\in(R[z]/(z^{t}%
))[x]
\]
of $B$ determined by the following recursion: $c_{0}=1$ and%
\[
c_{k}=-\frac{1}{k}\left(  c_{k-1}\mathrm{tr}(B)+c_{k-2}\mathrm{tr}%
(B^{2})+\cdots+c_{1}\mathrm{tr}(B^{k-1})+c_{0}\mathrm{tr}(B^{k})\right)
\]
for $1\leq k\leq tn$ (Newton formulae). In view of%
\[
\mathrm{tr}(B^{k})=\mathrm{tr}((\mu_{n}(A))^{k})=\mathrm{tr}(\mu_{n}%
(A^{k}))\in R^{\sigma},
\]
we deduce that $c_{i}\in R^{\sigma}$ for each $0\leq i\leq tn$. Thus
$\det(xI-B)\in R^{\sigma}[x]$ and the Cayley-Hamilton identity for $B$ is of
the form%
\[
B^{tn}+c_{1}B^{tn-1}+\cdots+c_{tn-1}B+c_{tn}I_{n}=0.
\]
Notice that for $r\in R$ and $c\in R^{\sigma}$ we have $c\sigma^{i-1}%
(r)=\sigma^{i-1}(cr)$ and $c\mu_{n}(A^{k})=\mu_{n}(cA^{k})$ follows from%
\[
\mu(cr_{0}+cr_{1}w+\cdots+cr_{t-1}w^{t-1}+(w^{t}))=c\mu(r_{0}+r_{1}%
w+\cdots+r_{t-1}w^{t-1}+(w^{t})).
\]
Thus%
\[
\left(  \mu_{n}(A)\right)  ^{tn}+c_{1}\left(  \mu_{n}(A)\right)
^{tn-1}+\cdots+c_{tn-1}\mu_{n}(A)+c_{tn}I_{n}=
\]%
\[
\mu_{n}(A^{tn}+c_{1}A^{tn-1}+\cdots+c_{tn-1}A+c_{tn}I_{n})=0
\]
holds in $\mathrm{M}_{tn}(R[z]/(z^{t}))$ and the injectivity of $\mu_{n}$
gives the desired identity. $\square$

\bigskip

\noindent3. SUPERMATRIX ALGEBRAS DETERMINED\ BY\ INVOLUTIONS

\bigskip

For an arbitrary endomorphism $\sigma:R\longrightarrow R$\ and for the integes
$1\leq k\leq n$ a matrix $A\in\mathrm{M}_{n}(R)$ is called a $(\sigma
,n,k)$-supermatrix if $A$ is of the shape%
\[
A=\left[
\begin{array}
[c]{cc}%
A_{1,1} & A_{1,2}\\
A_{2,1} & A_{2,2}%
\end{array}
\right]  ,
\]
where $A_{1,1}$ is a $k\times k$ and $A_{2,2}$ is an $(n-k)\times(n-k)$ square
block, while $A_{1,2}$ is a $k\times(n-k)$ and $A_{2,1}$ is an $(n-k)\times k$
rectangular block such that $\sigma(u)=u$ for each entry $u$\ of $A_{1,1}$ and
$A_{2,2}$ and $\sigma(u)=-u$ for each entry $u$\ of $A_{1,2}$ and $A_{2,1}$.
Let $\mathrm{M}_{n,k}^{\sigma}(R)$\ denote the set of $(\sigma,n,k)$%
-supermatrices. If $R=E$ is the (infinitely generated) Grassmann algebra and
$\tau(g_{0}+g_{1})=g_{0}-g_{1}$ is the natural $E\longrightarrow E$
involution, then $\mathrm{M}_{n,k}^{\tau}(E)$ is the classical algebra of
$(n,k)$-supermatrices (see [K]).

\bigskip

\noindent\textbf{3.1.Proposition.}\textit{ The set }$\mathrm{M}_{n,k}^{\sigma
}(R)$\textit{ is a subring (subalgebra) of }$\mathrm{M}_{n}(R)$\textit{.}

\bigskip

\noindent\textbf{Proof.} Straightforward verification. $\square$

\bigskip

\noindent\textbf{3.2.Theorem.}\textit{ Let }$\frac{1}{2}\in R$\textit{ and
}$\sigma:R\longrightarrow R$\textit{ be an arbitrary endomorphism. For }$r\in
R$\textit{ the definition}%
\[
\Theta(r)=\frac{1}{2}\left[
\begin{array}
[c]{cc}%
r+\sigma(r) & r-\sigma(r)\\
r-\sigma(r) & r+\sigma(r)
\end{array}
\right]
\]
\textit{gives an embedding }$\Theta:R\longrightarrow\mathrm{M}_{2}%
(R)$\textit{. If }$\sigma$\textit{\ is an involution (}$\sigma^{2}%
=1$\textit{), then }$\Theta$\textit{ is an }$R\longrightarrow\mathrm{M}%
_{2,1}^{\sigma}(R)$\textit{\ embedding.}

\bigskip

\noindent\textbf{Proof.} We give the details of the straightforward proof.

\noindent The additive property of $\Theta$ is clear. In order to prove the
multiplicative property of $\Theta$ take $r,s\in R$ and compute the product of
the $2\times2$ matrices $\Theta(r)$ and $\Theta(s)$:%
\[
\Theta(r)\cdot\Theta(s)=\frac{1}{4}\left[
\begin{array}
[c]{cc}%
r+\sigma(r) & r-\sigma(r)\\
r-\sigma(r) & r+\sigma(r)
\end{array}
\right]  \cdot\left[
\begin{array}
[c]{cc}%
s+\sigma(s) & s-\sigma(s)\\
s-\sigma(s) & s+\sigma(s)
\end{array}
\right]  =
\]%
\[
\frac{1}{4}\!\left[
\begin{array}
[c]{cc}%
\!(r\!+\!\sigma(r))(s\!+\!\sigma(s))\!+\!(r\!-\!\sigma(r))(s\!-\!\sigma
(s))\! & \!(r\!+\!\sigma(r))(s\!-\!\sigma(s))\!+\!(r\!-\!\sigma
(r))(s\!+\!\sigma(s))\!\\
\!(r\!-\!\sigma(r))(s\!+\!\sigma(s))\!+\!(r\!+\!\sigma(r))(s\!-\!\sigma
(s))\! & \!(r\!-\!\sigma(r))(s\!-\!\sigma(s))\!+\!(r\!+\!\sigma
(r))(s\!+\!\sigma(s))\!
\end{array}
\right]  \!=
\]%
\[
\frac{1}{4}\left[
\begin{array}
[c]{cc}%
2rs+2\sigma(r)\sigma(s) & 2rs-2\sigma(r)\sigma(s)\\
2rs-2\sigma(r)\sigma(s) & 2rs+2\sigma(r)\sigma(s)
\end{array}
\right]  =\frac{1}{2}\left[
\begin{array}
[c]{cc}%
rs+\sigma(rs) & rs-\sigma(rs)\\
rs-\sigma(rs) & rs+\sigma(rs)
\end{array}
\right]  =\Theta(rs).
\]
The injectivity of $\Theta$ follows from the fact that%
\[
r+\sigma(r)=s+\sigma(s)\text{ and }r-\sigma(r)=s-\sigma(s)
\]
imply%
\[
2r=(r+\sigma(r))+(r-\sigma(r))=(s+\sigma(s))+(s-\sigma(s))=2s.
\]
If $\sigma$\ is an involution, then%
\[
\sigma(r+\sigma(r))=\sigma(r)+\sigma^{2}(r)=\sigma(r)+r\text{ and }%
\sigma(r-\sigma(r))=\sigma(r)-\sigma^{2}(r)=\sigma(r)-r
\]
ensure that $\Theta(r)\in\mathrm{M}_{2,1}^{\sigma}(R)$. $\square$

\bigskip

\noindent\textbf{3.3.Remark.} Taking $R=\mathrm{M}_{n}(E)$ and the natural
extension $\sigma=\tau_{n}$, we obtain the the well known embedding%
\[
\Theta:\mathrm{M}_{n}(E)\longrightarrow\mathrm{M}_{2,1}^{\tau_{n}}%
(\mathrm{M}_{n}(E))\cong\mathrm{M}_{2n,n}^{\tau}(E).
\]

\bigskip

\noindent REFERENCES

\bigskip

\noindent\lbrack Do] M. Domokos, \textit{Eulerian polynomial identities and
algebras satisfying a standard identity}, Journal of Algebra 169(3) (1994), 913-928.

\noindent\lbrack Dr] V. Drensky, \textit{Free Algebras and PI-Algebras,}
Springer-Verlag, 2000.

\noindent\lbrack DrF] V. Drensky and E. Formanek, \textit{Polynomial Identity
Rings,} Birkh\"{a}user-Verlag, 2004.

\noindent\lbrack K] A. R. Kemer,\textit{\ Ideals of Identities of Associative
Algebras,} Translations of Math. Monographs, Vol. 87 (1991), AMS, Providence,
Rhode Island.

\noindent\lbrack MMSzvW] L. M\'{a}rki, J. Meyer, J. Szigeti and L. van Wyk,
\textit{Matrix representations of finitely generated Grassmann algebras and
some consequences,} arXiv:1307.0292

\noindent\lbrack SSz] S. Sehgal and J. Szigeti, \textit{Matrices over
centrally }$\mathbb{Z}_{2}$\textit{-graded rings,} Beitrage zur Algebra und
Geometrie (Berlin) 43(2) (2002), 399-406.
\end{document}